\documentclass[10pt]{article}

\usepackage{amsmath, amsthm, amsfonts, amssymb}
\usepackage{mathrsfs}
\usepackage{amssymb,amsmath,latexsym}
\usepackage{setspace}
\onehalfspace

\newtheorem{definition}{Definition}[section]
\newtheorem{lemma}{Lemma}[section]
\newtheorem{theorem}{Theorem}[section]
\newtheorem{corollary}{Corollary}[section]
\newtheorem{ex}{Example}[section]
\newtheorem{proposition}{Proposition}[section]

\def\blemma{\begin{lemma}\sl{}\def\elemma{\end{lemma}}}
\def\btheorem{\begin{theorem}\sl{}\def\etheorem{\end{theorem}}}

\def\bdefinition{\begin{definition}\sl{}\def\edefinition{\end{definition}}}

\def\be{\begin{equation}}\def\ee{\end{equation}}
\def\bee{\begin{equation*}}\def\eee{\end{equation*}}
\def\beqlb{\begin{eqnarray}}\def\eeqlb{\end{eqnarray}}
\def\beqnn{\begin{eqnarray*}}\def\eeqnn{\end{eqnarray*}}

\def\bp{\begin{proposition}}\def\ep{\end{proposition}}

\def\qed{\hfill$\Box$\medskip}

\def\<{\langle}\def\>{\rangle}

\def\E{\mathbb{E}}
\def\FF{\mathcal{F}}

\def\G{{\cal G}}
\def\H{{\cal H}}

\def\IR{\mathbb{R}}

\def\IR{{I\!\!R}}

\def\IN{{I\!\!N}}

\begin{document}
\title{\noindent{\Large\bf Weak Convergence of Equity Derivatives Pricing with Default Risk \footnote{
This work is partially supported by NSFC grants (No.71073129, No.11426187, No.11126236 and No.11201150),  China Postdoctoral Science Foundation (No. 2013M541526), JBK.120211, the Fundamental Research Funds for Central Universities
(No.239201278210075), and the 111 Project(B14019).}}
\author{\ Gaoxiu Qiao \footnote{The School of Finance, Southwestern University of Finance and
Economics, Chengdu, Sichuan Province, 611130, P. R. China. Email:
qiaogaoxiu@163.com}\ \ \ \
 and ~ Qiang Yao \footnote{Corresponding Author. School of Finance and Statistics, East China Normal University, Shanghai 200241,
 P. R. China. Email: qyao@sfs.ecnu.edu.cn}
 }}
\date{}
\maketitle

\noindent{\bf{Abstract}}\ This paper presents a discrete--time equity
derivatives pricing model with default risk in a no--arbitrage
framework. Using the equity--credit reduced form approach where
default intensity mainly depends on the firm's equity value, we
deduce the Arrow--Debreu state prices and the explicit pricing result
in discrete time after embedding default risk in the pricing model.
We prove that the discrete--time defaultable equity derivatives
pricing has convergence stability, and it converges weakly to the
continuous--time pricing results.



\section{Introduction}

\setcounter{equation}{0}

Default risk is the risk that the agents cannot fulfill their
obligations in the contracts. The reduced form approach has become a standard tool for modeling default risk. It considers the default to be an exogenously
specified jump process, and derives the default probability as the
instantaneous likelihood of default, see, for example, Jarrow and Turnbull \cite{jarrow2012pricing}, Duffie and Singleton \cite{duffie1999modeling}, Lando
\cite{lando1998cox}. The default time is usually defined as the first
jump time of a Cox process with a given intensity~(hazard rate).
Hence, these models are frequently called \emph{intensity models}.

Recently, an alternative model named \emph{equity--credit market
approach} has emerged. It assumes that the default intensity depends on the
firm's equity value(stock prices) and allows the stock price to jump
to zero at the time of default. It has both reduced form and
structural features. Default risk is incorporated in this equity modeling
approach by assuming that the stock price $S_t$ at time $t$ can jump
to zero with an intensity, which is assumed to be a function of
$S_t$. The models described above are all continuous--time models, they are widely used to model default risk.

However, continuous--time models are often too complicated to handle,
it is necessary to deduce discrete--time models and show that the
pricing processes converge to the continuous--time models. This is not a trivial job, since weak convergence, by its nature, is not tied to a single probability space. Some authors have presented different discrete--time models for
derivatives pricing and have established some weak convergence results. See, for example, Cox, Ross and Rubinstein \cite{cox1979option}, He \cite{he1990convergence}, Duffie and Protter \cite{duffie2006discrete}, Nieuwenhuis and
Vellekoop \cite{nieuwenhuis2004weak}, etc.

In this paper, our aim is to present a discrete--time equity
derivatives pricing model with default risk in a no--arbitrage
framework, and prove that the pricing in discrete--time converges
weakly to the continuous--time pricing results. In comparison, our
method is different from Nieuwenhuis and Vellekoop \cite{nieuwenhuis2004weak}.
Following the discrete framework of He \cite{he1990convergence} and equity--credit
market approach presented in \cite{bielecki2008defaultable}, we describe the
discrete--time pricing model in a no--arbitrage framework. After
embedding default risk, we deduce the Arrow--Debreu state prices and
the explicit pricing result in discrete time. In order to prove the
weak convergence of pricing processes, several auxiliary results are
presented.

The paper is organized as follows: In section 2, we introduce the
continuous--time model using equity--credit reduce-from approach; In
section 3, we illustrate a discrete--time model of the equity
derivatives pricing with default risk; In section 4, weak
convergence of equity derivatives pricing with default risk from
discrete--time to continuous--time pricing is proved; Finally, in section
5, we summarize the article and make concluding remarks.

\section{The continuous--time model}

\setcounter{equation}{0}

We first recall the continuous--time defaultable contingent claims
pricing model. Given a probability space $(\Omega,\FF,P)$,
$T$ is a strictly positive real number which represents the final date, $(\omega_t)_{0\leq t\leq T}$
is a Brownian motion. Let $\FF_{t}=\sigma(\omega_s,s\leq t)$ for $t\geq0$. We suppose
$\FF_{t}\subset \FF$ for all $t$, and $P$ is the real-world
probability. Furthermore, we denote by ``$\Rightarrow$'' weak convergence from now on.

A default event occurs at a random time $\tau$, where $\tau$ is a non-negative
random variable. The default process is defined as
$N_t\triangleq\textbf{1}_{\{\tau\leq t\}}$, and $\H_t=\sigma(N_s, s\leq t)$,
the filtration $\H$ is used to describe the information about
default time, where $\H=\bigcup\limits_{0\leq t\leq T}\H_t$. At any time $t$,
the agent's information on the securities prices and default time is
$\G_t=\FF_t\vee\H_t$ and the agent knows whether or not the default
has appeared. Hence, the default time $\tau$ is a $\G$ stopping time
where $\G=\bigcup\limits_{0\leq t\leq T}\G_t$. In fact, $\G$ is the smallest
filtration which contains $\FF$ and allows $\tau$ to be a stopping
time. Assume that the pre--default stock price $S_t$ has the
following dynamics \beqlb\label{e:2.1}
dS_t=(b(S_t)+\lambda(S_t,t)S_t)dt+\sigma(S_t,t)S_td\omega_t,\quad
S_0>0. \eeqlb Here we assume that $b(x)$ is continuous,
$\sigma(S,t)$ is a positively bounded and nonsingular Borel-measurable function. In particular we have that $\sigma(S,t)\geq \sigma$
for some positive constant $\sigma$, $\lambda(S,t)$ is a nonnegative, bounded,
continuous, $\FF$--progressively measurable and integrable function. The
functions $b(S)$, $\lambda(S,t)S$ and $\sigma(S,t)S$ are Lipschitz
continuous in $S$, uniformly in $t$.

The bond price $B_t$ satisfies $dB_t=B_tr(S_t)dt$ and $B_0=1$,
where $r(x)$ is a nonnegative continuous function, representing the
riskyless interest rate. Suppose there exists a constant $K>0$ such that
$|x^2r(x)|\leq K(1+x^2)$.

There exists a $\G$ equivalent martingale measure $Q^*$ which is
defined as $dQ^*|_{\FF_t}=\xi_tdP|_{\FF_t}$, where $\xi_t$ is the
Radon-Nikod$\acute{y}$m density satisfying
 \beqlb\label{e:2.2}
d\xi_t=\xi_{t}\theta(S_t)d\omega_t,\quad \xi_0=1. \eeqlb Here
$\theta(x)=-\sigma(x)^{-1}(b(x)-r(x)x)$. Define $W_t$ via
$dW_t=d\omega_t-\theta(S_t)dt$, then $W_t$ is a Brownian motion with
respect to $\FF$, and under the changed measure \beqlb\label{e:2.3}
dS_t=S_t[(r(S_t)+\lambda(S_t,t))dt+\sigma(S_t,t)dW_t],\quad S_0>0.
\eeqlb

Define $G_t\triangleq Q^*(\tau>t~|~\FF_t)$, $\Gamma_t\triangleq-\ln
G_t$. We call $\Gamma_t$ the \emph{$\FF$ hazard process of $\tau$}. For
the detailed properties, one can refer to Bielecki and Rutkowski
\cite{bielecki2004credit}.

Let $g(\cdot):\IR\rightarrow \IR$ be a square integrable and
measurable function, the equity derivatives are defined to be
securities that pay $g(S_T)$ dollars on the final date. This
formulation subsumes all of the usual examples, such as the European
options, convertible bonds and so on. The prices of equity
derivatives at time $t$ are \beqlb\label{e:2.4}
V(S_t,t)=\textbf{1}_{\{\tau>t\}}\E_{Q^*}\left[\left.\frac{B_te^{\Gamma_t}}{B_Te^{\Gamma_T}}g(S_T)~\right|~\FF_t\right].\eeqlb

Poisson process with stochastic intensity is called \emph{Cox
process}. Given $\lambda(S_u,u)$, denote by $\left\{\overline{C}_t\right\}$ the Poisson process with intensity $C_t=\int_0^t
\lambda(S_u,u)du$. Then $\left\{\overline{C}_t\right\}$ is a Cox process. Following the equity--credit market
models, the canonical construction of default time $\tau$ under the Cox
process $\left\{\overline{C}_t\right\}$ is defined as
$\tau=\inf\{t\geq0:~C_t\geq\Theta\}$, where $\Theta\sim Exp(1)$ and is independent of $\FF$
under $Q^*$. Then
$$
Q^*(\tau>t~|~\FF_t)=Q^*(\Theta>C_t~|~\FF_t)=e^{-C_t}. $$

It is easy to see that under this condition, the default time is the
first jump time of the Cox process, the $\FF$ hazard process of $\tau$
satisfies $$\Gamma_t=-\ln Q^*(\tau>t~|~\FF_t)=-\ln
Q^*(\Theta>C_t~|~\FF_t)=C_t.$$

Let $\Delta$ denote the bankruptcy state when the firm defaults at
time $\tau$. Then we can also write the dynamics for the stock price
subject to bankruptcy $S_t^\Delta$ as follows:
$$
dS_t^\Delta =S_t^\Delta[r(S_t)dt + \sigma(S_t,t)dW_t-dM_t],
$$
where $M_t=N_t-\int_0^{t\wedge\tau}\lambda(S_u,u)du$, and $M_t$ is a
martingale. Moreover, referred to Hypothesis (H) in
Blanchet-Scalliet and Jeanblanc \cite{blanchet2004hazard}: all $\FF$-martingales are
$\G$-martingales. It implies that the $\FF$-Brownian motion $W_t$ remains
a Brownian motion under the extended probability measure $Q^*$ and
with respect to the enlarged filtration $\G$ and is independent of
$M_t$.

Then we have \beqnn
V(S_t,t)=\textbf{1}_{\{\tau>t\}}\E_{Q^*}\left[\left.e^{-\int_t^T(r_s+\lambda_s)ds}g(S_T)~\right|~\FF_t\right].
\eeqnn Here we write $\sigma_t=\sigma(S_t)$, $r_t=r(S_t)$,
$\lambda_t=\lambda(S_t,t)$ for simplicity. We can obtain the
following result.
 \blemma\label{l:2.1}
Let
$Y(S_t,t)=\E_{Q^*}\left[\left.e^{-\int_t^T(r_s+\lambda_s)ds}g(S_T)~\right|~\FF_t\right]$.
Then it satisfies
 \beqlb\label{e:2.5}
 \frac{\partial Y}{\partial
t}+\frac{\sigma_t^2S_t^2}{2}\frac{\partial^2Y}{\partial
{S^2}}+(r_t+\lambda_t)S_t\frac{\partial Y}{\partial
S}-(r_t+\lambda_t)Y = 0. \eeqlb \elemma \noindent{Proof.} Let
$\hat{Y}(S_t,t) = E_{Q^*}\left[\left.\frac{G_Tg(S_T)}{B_T}~\right|~\FF_t\right]$.
It can be regarded as the discount price of contingent claim
$G_Tg(S_T)$ at time $t$, then it is $\FF$ martingale. By $It\hat{o}$'s formula,
$$
\frac{\partial\hat{Y}}{\partial
t}+\frac{\sigma_t^2S_t^2}{2}\frac{\partial^2\hat{Y}}{\partial
S^2}+(r_t+\lambda_t)S_t\frac{\partial\hat{Y}}{\partial S}=0.$$ Since
 $\hat{Y}(S_t,t)=e^{-\int_0^t(r_s+\lambda_s)ds}Y(S_t,t)$, (\ref{e:2.5}) is proved after using $It\hat{o}$'s formula again.\qed


\section{Discrete--time model in the defaultable market}
 \setcounter{equation}{0}

For simplicity, the time horizon is assumed to be $[0,1]$ and
divided into $n$ steps, the length of every step is $1/n$. For $k=1,2,\cdots,n$, let
$\varepsilon^k$ be a random variable on the probability space
$\bar{\Omega}=\{\omega_1,\omega_2\}$. For example, set
$\varepsilon^k(\omega_1)=1$,  $\varepsilon^k(\omega_2)=-1$,  and
$P(\{\omega_1\})=P(\{\omega_2\})=\frac{1}{2}$.

Let $\bar{\Omega}_n=\overbrace{\bar{\Omega}\times\bar{\Omega}
\times\cdot\cdot\cdot\times\bar{\Omega}}^n=
\{\omega_1,\omega_2\}^n$, $P_n=\overbrace{P\times P \times ...\times
P}^n$. Then $P_n$ is the probability measures defined on
$\bar{\Omega}_n$, representing the real-world probability, $\FF_n$ is
the filtration generated by $\varepsilon^k, k=1,\ldots, n$,
$\{\varepsilon^1,\varepsilon^2,\ldots,\varepsilon^n\}$ is a
sequence of independent and identically distributed random vectors
defined on $\{\bar{\Omega}_n,\FF_n, P_n\}$.

There are two financial assets in the market: stock and bond. Since
the increment of Brownian motion can be approximated by a sequence
of independent and identically distributed random variables, the
pre--default stock prices and bond prices can be written as \beqnn
&&S_{k+1}^n=S^k_{n}+\frac{b(S_k^n)+\lambda(S_k^n,\frac{k}{n})S^k_{n}}{n}+\frac{\sigma(S_k^n,\frac{k}{n})S^k_{n}}{\sqrt{n}}\varepsilon^{k+1},\quad
S_0^n=S_0.\\
&&B_{k+1}^n=B_k^n\left(1+\frac{r(S_k^n)}{n}\right),\quad B_0^n=1.
  \eeqnn
Here $S_k^n$, $B_k^n$ denote the stock prices and bond prices at
time $\frac{k}{n}$ respectively. Let $\tilde{S}_t^n=S_{[nt]}^n$,
$\tilde{B}_t^n=B_{[nt]}^n$. Then $\tilde{S}_t^n$,  $\tilde{B}_t^n$
are Markov processes and have jumps only at time $\frac{k}{n}$.
Moreover $\tilde{S}_t^n$ can be expressed as follows
 \beqnn
 \tilde{S}_t^n&=&S_0+\sum_{i=0}^{[nt]-1}\frac{b(S_i^n)+\lambda(S_i^n,\frac{i}{n})S_i^n}{n}
 +\sum_{i=0}^{[nt]-1}\frac{\sigma(S_i^n, \frac{i}{n})S_i^n}{\sqrt{n}}\varepsilon^{i+1}\\
 &=&S_0+\int_0^{\frac{[nt]}{n}}(b(\tilde{S}_u^n,u)+\lambda(\tilde{S}_u^n,u)\tilde{S}_u^n)du+
 \sum_{i=0}^{[nt]-1}\frac{\sigma(S_i^n,\frac{i}{n})S_i^n}{\sqrt{n}}\varepsilon^{i+1}.
 \eeqnn
 Similarly,
$$\tilde{B}_t^n=B_0+\sum_{i=0}^{[nt]-1}\frac{r(S_i^n)B_i^n}{n}=B_0+\int_0^{\frac{[nt]}{n}}r(\tilde{S}_u^n)\tilde{B}_u^ndu.$$

The above discrete framework is employed in He \cite{he1990convergence} where he
assumes the stock prices satisfying
$\pi(\omega_1;S_k^n)S_{k+1}^n(\omega_1)
+\pi(\omega_2;S_k^n)S_{k+1}^n(\omega_2)=S_k^n$, where
$\pi(\omega_s;S_k^n)(s=1, 2)$ are considered as the Arrow--Debreu
state prices and the discount stock prices are martingale. He
\cite{he1990convergence} gives the result that there exists unique equivalent
martingale measure $Q_n$ in discrete--time defaultable market, and
$dQ_n=\xi_n^ndP_n$, where $\xi_k^n=2^{k}\pi_k^n B_k^n$, $\pi_k^n$ is
defined as the product of Arrow--Debreu state prices from 0 to k, the
default-free discrete--time market is complete. Let
$\tilde{\xi}_t^n=\xi_{[nt]}^n$, then $\displaystyle{\tilde{\xi}_t^n=\xi_0+\sum_{i=0}^{[nt]-1}\frac{\theta(S_i^n)\xi_i^n}{\sqrt{n}}\varepsilon^{i+1}}$.

Suppose  default occurs at random time $\tau_n$, where $\tau_n$ is a non-negative
random variable. The default process is defined as
$N_k^n=1_{\{\tau_n\leq \frac{k}{n}\}}$,
 $\sigma$ filtration $\H_k^n=\sigma(N_i^n, 0\leq i\leq k)$ and $\H^n$ is used to
describe the information about default time. At time $\frac{k}{n}$,
the agent's information on the prices and on the default time is
$\G_k^n=\FF_k^n\vee\H_k^n$. Hence, the default time $\tau_n$ is a
$\G^n$ stopping time where $\G^n=\{\G_k^n, 0\leq k\leq n\}$.

From Blanchet-Scalliet and Jeanblanc \cite{blanchet2004hazard}, if the defaultfree market
is complete and arbitrage-free, the defaultable market is
arbitrage-free, then there exists an equivalent martingale measure
$Q_n^*$ in $\G^n$-market.

\bdefinition\label{d:3.1} For any $0\leq k \leq n$, Let $Q_n^*(\tau_n=0)=0$,
$Q_n^*(\tau_n>\frac{k}{n})>0$. We write $F_k^n=Q_n^*(\tau_n\leq
\frac{k}{n}~|~\FF_k^n)$,
$G_k^n=1-F_k^n=Q_n^*(\tau_n>\frac{k}{n}~|~\FF_k^n)$. Suppose $F_k^n<1$,
Let $\Gamma_k^n\triangleq-\ln G_k^n=-\ln(1-F_k^n)$ is called the
$\FF^n$ hazard process of $\tau_n$ under $Q_n^*$. \edefinition

Several properties in the continuous--time model still hold here, such as
$(F_k^n)$ is nonnegative bounded submartingale, $\Gamma_k^n$ is increasing,
 and $L_k^n=1_{\{\tau_n>\frac{k}{n}\}}e^{\Gamma_k^n}$ is a martingale,
 the detailed properties can refer to \cite{blanchet2004hazard}.

Assume the Cox process is defined via intensity $\tilde{C}_t\triangleq\int_0^t
\lambda(\tilde{S}_u^n,u)du$, $\Theta$ is the random variable with an
exponential law of parameter 1, which is independent of $\FF^n$ under
$Q_n^*$, then we can give the canonical construction of $\tau_n$,
$\tau_n:\bar{\Omega}_n\rightarrow [0, T],
\tau_n=\inf\left\{\frac{[nt]}{n}\geq0;
\tilde{C}_{\frac{[nt]}{n}}\geq\Theta\right\}$. Therefore,
$$
Q_n^*\left(\left.\tau_n>\frac{[nt]}{n}~\right|~\FF_{[nt]}^n\right)
=Q_n^*\left(\left.\Theta>\tilde{C}_{\frac{[nt]}{n}}~\right|~\FF_{[nt]}^n\right)
=e^{-\tilde{C}_{\frac{[nt]}{n}}}.
$$
The $\FF^n$ hazard process $\tau_n$ satisfy
$$
\Gamma_k^n=-\ln Q_n^*\left(\left.\tau_n>\frac{k}{n}~\right|~\FF_k^n\right)=-\ln
Q_n^*\left(\left.\Theta>\tilde{C}_{\frac{k}{n}}~\right|~\FF_{k}^n\right)=\tilde{C}_{\frac{k}{n}}.$$
Moreover,
$$\Gamma_{k+1}^n=\Gamma_k^n+\int_{\frac{k}{n}}^{\frac{k+1}{n}}
\lambda(\tilde{S}_u^n,u)du=\Gamma_k^n+\frac{\lambda(S_k^n,\frac{k}{n})}{n}.
$$ Define $\tilde{\Gamma}_t^n=\Gamma_{[nt]}^n$, then
$\tilde{\Gamma}_t^n$ is a sequence of Markov process on probability
space$(\bar{\Omega}_n, \FF^n, Q_n^*)$ with sample path in $D_{\IR}[0,
1]$ and $$
 \tilde{\Gamma}_t^n=\sum_{i=0}^{[nt]-1}\frac{\lambda(S_i^n,\frac{i}{n})}{n}=\int_0^{\frac{[nt]}{n}}\lambda(\tilde{S}_u^n,u)du,~~~~\tilde{\Gamma}_0^n=0.$$ Although $\tau_n$ and $\tau$ are defined in
different probability spaces, the canonical construction provides us
a feasible way to prove the weak convergence of the default process
which will be shown in the next section.

Define auxiliary discount process $$\beta_k^n=B_k^ne^{\Gamma_k^n}\triangleq\prod_{i=0}^{k-1}\left(1+\frac{\tilde{r}_i}{n}\right),$$ then we have $\displaystyle{\frac{\beta_{k+1}^n}{\beta_k^n}=\frac{B_{k+1}^ne^{\Gamma_{k+1}^n}}{B_k^ne^{\Gamma_k^n}}
=\left(1+\frac{r(S_k^n)}{n}\right)\exp\left(\frac{\lambda(S_k^n)}{n}\right)}$,
where  $\displaystyle{\tilde{r}_k= (r_k + \lambda_k) + \frac{\lambda_k^2 +
2r_k\lambda_k}{2n}+o\left(\frac{1}{n}\right)}$. Here we write $r_k=r(S_k^n)$ and
$\lambda_k=\lambda(S_k^n)$ for simplicity.
Suppose the equity prices and bond prices satisfy $$
\tilde{\pi}(\omega_1;S_k^n)S_{k+1}^n(\omega_1)
+\tilde{\pi}(\omega_2;S_k^n)S_{k+1}^n(\omega_2)=S_k^n,~~~\tilde{\pi}(\omega_1;S_k^n)B_{k+1}^n(\omega_1)
+\tilde{\pi}(\omega_2;S_k^n)B_{k+1}^n(\omega_2)=B_k^n.$$ Solve
the above two equations, we get $\displaystyle{\tilde{\pi}(\omega_s;S_k^n)=\frac{1}{2}  \left(1+\frac{\theta(S_k^n)}{\sqrt{n}}
\varepsilon^{k+1}\right)\left(1+\frac{\tilde{r}_k}{n}\right)^{-1}}$.
Set $$\tilde{\pi}_k^n=\tilde{\pi}(\cdot;S_{k-1}^n)\tilde{\pi}
(\cdot;S_{k-2}^n)\cdots\tilde{\pi(}\cdot;S_0^n),~~k=1,2,\ldots,n~~~~~~\tilde{\pi}_0^n=1.$$ Since $\theta$ is bounded, for $n$ large
enough, $\tilde{\pi}$ is non-negative. Then
$$\tilde{\pi}(\omega_1;S_k^n)
+\tilde{\pi}(\omega_2;S_{k}^{n})=\left(1+\frac{\tilde{r}_k}{n}\right)^{-1},~~~\xi_k^n=2^n\tilde{\pi}_k^n\beta_k^n=2^n\pi_k^nB_k^n.$$

The price of defaultable contingent claims $g(S_n^n)$ at time
$\frac{k}{n}$ is \beqlb\label{e:3.1}
V^n\left(S_k^n,\frac{k}{n}\right)=1_{\{\tau_n>\frac{k}{n}\}}\E_{Q_n}\left[\left.\frac{B_k^ne^{\Gamma_k^n}}{B_n^ne^{\Gamma_n^n}}g(S_n^n)~\right|~\FF_k^n\right].
\eeqlb
 \blemma\label{l:3.1} Let $Y_n\left(S^n_k,\frac{k}{n}\right)=\E_{Q_n}
\left[\left.\frac{\beta_k^n}{\beta_n^n}g(S_n^n)~\right|~\FF_k^n\right]$. Then the
following equation holds:
 \beqlb\label{e:3.2}
 Y_n\left(S_k^n,\frac{k}{n}\right)&=\tilde{\pi}(\omega_1;S_k^n)Y_n\left(S_{k+1}^n(\omega_1),\frac{k}{n}\right)+
\tilde{\pi}(\omega_2;S_k^n)Y_n\left(S_{k+1}^n(\omega_2),\frac{k}{n}\right).\eeqlb
\elemma \noindent{Proof.} Since $Y_n\left(S_k^n,\frac{k}{n}\right)$ is
a $\FF^n$ martingale, then we have \begin{align*}
&Y_n\left(S_k^n,\frac{k}{n}\right)=\E_{Q_n}\left[\left.Y_n\left(S_{k+1}^n,\frac{k+1}{n}\right)\left(1+\frac{\tilde{r}_k}{n}\right)^{-1}~\right|~\FF_k^n\right]\\
=&\E_{P_n}\left[\left.\frac{\xi_{k+1}^n}{\xi_k^n}Y_n\left(S_{k+1}^n,\frac{k+1}{n}\right)\left(1+\frac{\tilde{r}_k}{n}\right)^{-1}~\right|~\FF_k^n\right]=\E_{P_n}\left[\left.\frac{2\pi_{k+1}^nB_{k+1}^n}{\pi_k^nB_k^n}Y_n\left(S_{k+1}^n,\frac{k+1}{n}\right)\left(1+\frac{\tilde{r}_k}{n}\right)^{-1}~\right|~\FF_k^n\right]\\
=&\frac{1}{2}\left(1+\frac{\theta(S_k^n)}{\sqrt{n}}\right)\left(1+\frac{\tilde{r}_k}{n}\right)^{-1}
Y_n\left(S_{k+1}^n(\omega_1),\frac{k+1}{n}\right)+\frac{1}{2}\left(1-\frac{\theta(S_k^n)}{\sqrt{n}}\right)
\left(1+\frac{\tilde{r}_k}{n}\right)^{-1}Y_n\left(S_{k+1}^n(\omega_2),\frac{k+1}{n}\right)\\
=&\tilde{\pi}(\omega_1;S_k^n)Y_n\left(S_{k+1}^n(\omega_1),\frac{k+1}{n}\right)
+\tilde{\pi}(\omega_2;S_k^n)Y_n\left(S_{k+1}^n(\omega_2),\frac{k+1}{n}\right).
 \end{align*}
Therefore, (\ref{e:3.2}) is proved.\qed

We can conclude that $\tilde{\pi}(\cdot;S_k^n)$ can be regarded as the
discrete--time Arrow--Debreu state prices in the defaultable market.

\section{The main result: Weak Convergence}
\setcounter{equation}{0}

In this section, we will prove the weak convergence of pricing
process for defaultable equity derivatives under the above model.
Firstly, we introduce infinite dimensional multiplicative
probability space
$\bar{\Omega}_\IN\triangleq\overbrace{\bar{\Omega}\times\bar{\Omega}
\times\cdot\cdot\cdot\times\bar{\Omega}}^\infty$, then
 $\bar{\Omega}_n$ is a subspace of $\bar{\Omega}_\IN$. From the
infinite multiply probability existence theorem, there exists a
unique probability measure $\hat{P}$ satisfying condition:
$\hat{P}(A\times\bar{\Omega}_{\IN \setminus I_n})=Q_n(A)$, where
$A\in \bar{\Omega}_n$, $I_n=\{1, 2, \ldots, n\}$.

He \cite{he1990convergence} proves the weak convergence of  Markov process vector
including equity prices, bond prices, Radon-Nikod$\acute{y}$m
density. Now we extend this result in the defaultable market.
Combining Martingale central limit theorem developed by Ethier and
Kurtz \cite{ethier1986markov} (Page 354), we get a similar result. \blemma\label{l:4.1}  For any $t\in[0,
1]$, let $\tilde{Z}_t^n=(\tilde{S}_t^n, \tilde{B}_t^n,
\tilde{\xi}_t^n, \tilde{\Gamma}_t^n )$, $Z_t=(S_t, B_t, \xi_t,
\Gamma_t)$. Then
 $$\tilde{Z}_{\cdot}^n\Rightarrow Z_{\cdot}~\text{as}~n\rightarrow\infty.$$
\elemma

Recall the definition of $\tau_n$ and $\tau$, the following
conclusion holds.

\blemma\label{l:4.2} For $t\in[0,1]$, let $X_n(t)=1_{\{\tau_n>\frac{[nt]}{n}\}}$,
$X(t)=1_{\{\tau>t\}}$. Then $X_n(\cdot)\Rightarrow X(\cdot)$.
\elemma \noindent{Proof.} For any $t\in[0,1]$, we have $\displaystyle{\sup_n
\E_{\hat{P}}\left[|X_n(t)|\right]=\sup_n \E_{Q_n}\left[|X_n(t)|\right]\leq
1<\infty}$. Therefore,
$$
\lim_{C\rightarrow\infty}\limsup\limits_n\hat{P}(|X_n(t)|>C)\leq\lim_{C\rightarrow\infty}\frac{\limsup\limits_n
\E_{\hat{P}}\left[|X_n(t)|\right]}{C}=0.$$ So $\{X_n(t)\}$is tight for any
$~t\in[0,1]$.

Choose sequences $\{\alpha_n\}$ and $\{\delta_n\}$ satisfying the following: for all $n$,
$\alpha_n$ is a stopping time with respect to the $\sigma$
filtration which is generated by the process $\{X_n(t): 0\leq
t\leq1\}$, and  $\alpha_n$ has only finite value; $\delta_n$ is a
constant and $0\leq\delta_n\leq1$. Moreover $\delta_n\rightarrow 0$,
as $n\rightarrow\infty$. \begin{align*}
&\hat{P}(|X_n(\alpha_n+\delta_n)-X_n(\alpha_n)|>\epsilon)=Q_n(|X_n(\alpha_n+\delta_n)-X_n(\alpha_n)|>\epsilon)\\
\leq&\frac{1}{\epsilon}\E_{Q_n}\left[|X_n(\alpha_n+\delta_n)-X_n(\alpha_n)|\right]\leq\frac{1}{\epsilon}\E_{Q_n}\left[1_{\{\tau_n>\frac{[n\alpha_n]}{n}\}}-1_{\{\tau_n>\frac{[n\alpha_n+n\delta_n]}{n}\}}\right]\\
=&\frac{1}{\epsilon}E_{Q_n}\left[\E_{Q_n}\left[\left.1_{\{\tau_n>\frac{[n\alpha_n]}{n}\}}~\right|~\FF_{[n\alpha_n]}^n\right]-
\E_{Q_n}\left[\left.1_{\{\tau_n>\frac{[n\alpha_n+n\delta_n]}{n}\}}~\right|~\FF_{[n\alpha_n+n\delta_n]}^n\right]\right]\\
=&\frac{1}{\epsilon}\E_{Q_n}\left[\exp(-\Gamma_{[n\alpha_n]}^n)-\exp(-\Gamma_{[n\alpha_n+n\delta_n]}^n)\right]\\
=&\frac{1}{\epsilon}\E_{Q_n}\left[\exp\left(-\sum_{i=0}^{[n\alpha_n]-1}\frac{\lambda(S_i^n)}{n}\right)\left(1-
\exp\left(-\sum_{i=[n\alpha_n]}^{[n\alpha_n+n\delta_n]-1}\frac{\lambda(S_i^n)}{n}\right)\right)\right].
\end{align*}
Since $\lambda(S,t)$ is a nonnegative \emph{bounded} continuous function, we have $\displaystyle{\sum_{i=[n\alpha_n]}^{[n\alpha_n+n\delta_n]-1}\frac{\lambda(S_i^n)}{n}\leq C\delta_n\rightarrow0}$ as $n\rightarrow\infty$, where $C$ is a fixed constant. Together with the fact that $\displaystyle{\exp\left(-\sum\limits_{i=0}^{[n\alpha_n]-1}\frac{\lambda(S_i^n)}{n}\right)\leq
1}$, we have $\displaystyle{\hat{P}(|X_n(\alpha_n+\delta_n)-X_n(\alpha_n)|>\epsilon)\rightarrow0}$,
that is, $
X_n(\alpha_n+\delta_n)-X_n(\alpha_n)\overset{\hat{P}}{\longrightarrow}0$
as $n\rightarrow\infty$.  By the criterion of Aldous \cite{aldous1978stopping} (Page 1),
$\{X_n(\cdot)\}$ is tight in $D_{\IR}[0,1]$.

We have $\displaystyle{X_n(t)=\textbf{1}_{\{\tau_n>\frac{[nt]}{n}\}}=\textbf{1}_{\{\Theta>\tilde{C}_{\frac{[nt]}{n}}\}}=\textbf{1}_{\{\Theta>\tilde{\Gamma}_t^n
\}}}$ and $\displaystyle{X(t)=\textbf{1}_{\{\tau>t\}}=\textbf{1}_{\{\Theta>C_t\}}=\textbf{1}_{\{\Theta>\Gamma_t\}}}$ for any $t\in[0,1]$.
Since $\tilde{\Gamma}^n\Rightarrow \Gamma$ as $n$ tends to infinity,
we have $\E_{Q_n}[e^{iuX_n(t)}]\rightarrow \E_{Q}[e^{iuX(t)}]$ as $n$ tends to infinity. According
to the dominated convergence theorem, for any
$t_1,t_2,\ldots,t_m\in[0,1]$, $u_1,\ldots, u_m\in \IR$,
$$
\E_{Q_n}\left[e^{i\sum_{j=1}^mu_jX_n(t_j)}\right]\rightarrow
\E_{Q}\left[e^{i\sum_{j=1}^m u_jX(t_j)}\right],\quad n\rightarrow \infty.
$$
Therefore, $\{X_n\}$ is tight, and their finite dimensional
distribution converges. From Ethier and Kurtz \cite{ethier1986markov} (Page 131),
$\{X_n\}$ converges weakly to $X$ .\qed

\blemma\label{l:4.3} For any integers $l,m, k\geq0$, $l\leq k\leq n$, there
exists a constant $C>0$, depending on $m$ (large enough), such that
$\E_{Q_n}\left[[S_k^n]^{2m}~|~\FF_l^n\right]\leq C(1+[S_0]^{2m})$. \elemma

 Next we prove
 the weak convergence of the second part in the equation of defaultable contingent claims
 prices, following a similar argument to the main theorem of
 He \cite{he1990convergence}.

\blemma\label{l:4.4} For any $t\in[0,1]$, let $$Y(t)\triangleq
Y(S_t,t)=\E_{Q}\left[\left.\frac{B_te^{\Gamma_t}}{B_Te^{\Gamma_T}}g(S_T)~\right|~\FF_t\right],~~~Y_n(t)\triangleq
Y_n\left(S^n_{[nt]},\frac{[nt]}{n}\right)=\E_{Q_n}\left[\left.\frac{B_{[nt]}^ne^{\Gamma_{[nt]}^n}}{B_n^ne^{\Gamma_n^n}}g(S_n^n)~\right|~\FF_{[nt]}^n\right].$$
 Suppose that $Y$ is continuously differentiable up to the third
order and that $Y$ and all of its derivatives up to the third order
satisfy a polynomial growth condition. Then $Y_n(\cdot)\Rightarrow
Y(\cdot)$ as $n$ tends to infinity. \elemma

{\bf Remark.} We get the idea of the proofs of Lemma \ref{l:4.3} and Lemma \ref{l:4.4} from He \cite{he1990convergence}, but the results in our paper are rather different
from them. In Lemma \ref{l:4.3}, we give the inequality for a more general
case. In Lemma \ref{l:4.4}, we prove that $\tilde{e}_t^n$ converges to zero
in the sense of almost everywhere.

Then combine Lemma \ref{l:4.2} and Lemma \ref{l:4.4}, we obtain the main result.

\btheorem\label{t:4.1} Suppose $g(\cdot)$ is $\IR\rightarrow \IR$ square
integrable measurable function. Let $\tilde{V}^n\left(\tilde{S}_t^n,
t\right)= V^n\left(S_{[nt]}^n,\frac{[nt]}{n}\right)$, $V(S_t,t)$ and
$\tilde{V}^n\left(\tilde{S}_t^n,\frac{[nt]}{n}\right)$ satisfy (\ref{e:2.4})
and (\ref{e:3.1}) respectively. Then we have
$$
\tilde{V}^n\left(\tilde{S}_{\cdot}^n, \cdot\right)\Rightarrow
V(S_{\cdot},\cdot) ,\quad as\quad n\rightarrow \infty.
$$
\etheorem

\noindent{Proof.} Clearly, we have $\displaystyle{\tilde{V}^n\left(\tilde{S}_t^n,\frac{[nt]}{n}\right)=X_n(t)Y_n(t)}$ and $\displaystyle{V(S_t,t)=X(t)Y(t)}$. By Lemma \ref{l:4.4}, $Y_n$ converges weakly to $Y$,
then $Y_n$ is relatively tight. Since $\IR$ is separable and $(\IR,
d)$ is complete, then $D_{\IR}[0,1]$ is separable, it follows that
$\{Y_n\}$ is tight.

By Lemma \ref{l:4.2}, $\{X_n\}$ is tight, together with the fact that
$Y(t)=Y(S_t,t)$ is continuous with respect to $~t$, then
$\{(X_n,Y_n)\}$ is tight according to Jacod-Shiryaev \cite{jacod1987limit} (Page 353).

Next, we only need  to prove the convergence of their finite
dimension distribution. That is, for any $~u_1,\ldots,u_m\in \IR$,
$~v_1,\ldots, v_m\in \IR,$
 \beqlb\label{e:4.1}
 \E_{Q_n}\left[e^{i\sum_{j=1}^m(u_jX_n(t_j)+v_jY_n(t_j))}\right]\rightarrow \E_{Q}\left[e^{i\sum_{j=1}^m
(u_jX(t_j)+v_jY(t_j))}\right],\quad n \rightarrow \infty. \eeqlb

From Lemma  4.2 and Lemma~4.4, we can obtain the convergence of the
finite dimension distribution of $\{X_n\}$, $\{Y_n\}$. Moreover,
$X_n(t)$ and $Y_n(t)$ are measurable with respect to $\H_t$ and
$\FF_t$ respectively, and $X_n(t)$, $Y_n(t)$ are independent, (\ref{e:4.1})
holds clearly. Then $(X_n,Y_n)\Rightarrow (X,Y)$ as
$n\rightarrow\infty$.

Let $f(X_n(t)Y_n(t))=X_n(t)Y_n(t)$, where $f$ is a continuous
function that maps $X_n(t)$, $Y_n(t)$ from $D_{\IR}[0,1]\times
D_{\IR}[0,1]$ to $D_{\IR}[0,1]$. By continuous mapping
theorem(Ethier-Kurtz \cite{ethier1986markov}, p.354), $X_n(\cdot)Y_n(\cdot)\Rightarrow
X(\cdot)Y(\cdot)$ as $n\rightarrow \infty$.\qed

In the following, we give the details of proofs of Lemma $4.1$,
$4.3$ and $4.4$. For simplicity, we write $b(\tilde{S}_u^n)=b_u$,
$\lambda(\tilde{S}_u^n,u)=\lambda_u$,
$\sigma(\tilde{S}_u^n)=\sigma_u$, $\theta(\tilde{S}_u^n)=\theta_u$
in the proof.

{\bf Proof of Lemma \ref{l:4.1}.} $d\Gamma_t=\lambda(S_t,t)dt$. Since
$b(\cdot)$, $\lambda(\cdot)S$, $\sigma(\cdot)S$ satisfy the
Lipschitz condition, therefore (\ref{e:2.1}) has a unique solution, which
implies that (\ref{e:2.2}), (\ref{e:2.3}) also have unique solutions
respectively. Since $\{\tilde{S}^n\}$, $\{\tilde{B}^n\}$,
$\{\tilde{\xi}^n\}$ and $\{\tilde{\Gamma}^n\}$ are processes that
are right continuous with left limits, hence $\{Z^n\}$ is a sequence
of Markov process vectors with sample path in $D_{\IR^4}[0,1]$,
where $D_{\IR^4}[0,1]$ is the space of functions from $[0,1]$ to
$\IR^4$, right continuous with left limits. Denote $L_t^n$ and
$A_t^n$ by \beqnn
L_t^n=\left(%
\begin{array}{c}
  \int_0^{\frac{[nt]}{n}}(b_u+\lambda_u\tilde{S}_u^n)du \\
  \int_0^{\frac{[nt]}{n}}r\tilde{B}_u^ndu\\
  0 \\
  \int_0^{\frac{[nt]}{n}}\lambda_u du \\
\end{array}%
\right),
A_t^n=\left(%
\begin{array}{cccc}
 \int_0^{\frac{[nt]}{n}}\sigma_u^2S^2du & 0
 & \int_0^{\frac{[nt]}{n}}\sigma_u S\theta_u\tilde{\xi}_u^ndu & 0 \\
  0 & 0 & 0 & 0 \\
   \int_0^{\frac{[nt]}{n}}\sigma_u S\theta_u\tilde{\xi}_u^ndu & 0&
    \int_0^{\frac{[nt]}{n}}(\theta_u\tilde{\xi}_u^n)^2du  & 0 \\
  0 & 0 & 0 & 0 \\
\end{array}%
\right) .\eeqnn Then $\{L^n\}$, $\{A^n\}$ are $4\times1$ and
$4\times4$(symmetric)matrix valued process respectively, and each of
their elements has a sample path in $D_{\IR}[0,1]$. Moreover,
$A_t^n-A_s^n$ is non-negative definite for $t>s\geq 0$. Define
$$\tau_n^q:=\inf\{t\leq T: |Z_t^n|\geq q~\text{or}~|Z_{t-}^n|\geq q \}.$$

Next, we prove the four conditions for martingale central limit
theorem holds.

(a) It is directly from $Z_0^n=Z_0=(S_0, B_0, \xi_0, 0)$.

(b) $
 M_t^n=\tilde{Z}^n-L_t^n=\left(%
\begin{array}{c}
  S_0+\sum_{i=0}^{[nt]-1}\frac{\sigma( S_i^n)S_i^n}{\sqrt{n}}\varepsilon^{i+1} \\
  B_0 \\
  \xi_0+\sum_{i=0}^{[nt]-1}\frac{\theta(S_i^n)\xi_i^n}{\sqrt{n}}\varepsilon^{i+1} \\
  0 \\
\end{array}%
\right).$

 Let $
N_k=S_0+\sum_{i=0}^{k-1}\frac{\sigma(S_i^n)S_i^n}{\sqrt{n}}\varepsilon^{i+1}$.
Clearly, $$\E_{Q_n}[N_{k+1}-N_k~|~S_k^n] =\E_{Q_n}\left[\left.\frac{\sigma(
S_k^n)S_k^n}{\sqrt{n}}\varepsilon^{k+1}~\right|~S_k^n\right]
=\frac{\sigma(
S_k^n)S_k^n}{\sqrt{n}}\E_{Q_n}[\varepsilon^{k+1}~|~S_k^n]=0.$$ So
$N_k$ are martingales. By the same arguments, $M_t^n$ and $M_t^n
(M_t^{n})^T-A_t^n$ are also martingales.

(c) $\displaystyle{\tilde{Z}_t^n-\tilde{Z}_{t-}^n=\left\{
\begin{array}{ll}
 \begin{pmatrix}\frac{b(S_{k-1}^n)+\lambda(S_{k-1}^n)S_{k-1}^n}{n}+\frac{\sigma(S_{k-1}^n)S_{k-1}^n}{\sqrt{n}}\varepsilon^k \\
   \frac{r(S_{k-1}^n)B_{k-1}^n}{n} \\
\frac{\theta(S_{k-1}^n)\xi_{k-1}^n}{\sqrt{n}}\varepsilon^k \\
  \frac{\lambda(S_{k-1}^n)}{n}\end{pmatrix}~~~~~~~~~~~~t=\frac{k}{n},\\~~~~~~~~~~~~~~~~~~~~~~~~~0~~~~~~~~~~~~~~~~~~~~~~~~~~~~~~~~~~~~~\frac{k}{n}<t<\frac{k+1}{n}.\end{array}\right.}$

By the definition of $\tau_n^q$, when $t\leq\tau_n^q$,
$|Z_t^n-Z_{t-}^n|\leq2q$, and $|Z_t^n-Z_{t-}^n|^2$ is of order
$\frac{1}{n}$. By the dominated convergence theorem, we have $\displaystyle{\lim_{n\rightarrow\infty}\E_n\left[\sup_{t\leq\tau_n^q}|Z_t^n-Z_{t-}^n|^2\right]=0}$.
Then using the same argument, we can get $\displaystyle{\lim_{n\rightarrow\infty}\E_n\left[\sup_{t\leq\tau_n^q}|L_t^n-L_{t-}^n|^2\right]=0}$ and $\displaystyle{\lim_{n\rightarrow\infty}\E_n\left[\sup_{t\leq\tau_n^q}|A_t^n-A_{t-}^n|\right]=0}$.

(d) For all $q>0$, we have
 \beqnn
L_t^n-\left(%
\begin{array}{c}
  \int_0^t(b(\tilde{S}_u^n)+\lambda(\tilde{S}_u^n)\tilde{S}_u^n)du \\
  \int_0^tr(\tilde{S}_u^n)\tilde{B}_u^ndu\\
  0 \\
  \int_0^t\lambda(\tilde{S}_u^n)du \\
\end{array}%
\right)=-\left(%
\begin{array}{c}
 (b(\tilde{S}_t^n)+\lambda(\tilde{S}_t^n)\tilde{S}_t^n)\left(t-\frac{[nt]}{n}\right) \\
 r(\tilde{S}_t^n)\tilde{B}_t^n\left(t-\frac{[nt]}{n}\right) \\
 0 \\
  \lambda(\tilde{S}_t^n)\left(t-\frac{[nt]}{n}\right) \\
\end{array}%
\right), \eeqnn where $\frac{k}{n}<t<\frac{k+1}{n}$, and
 the above equation equals to zero when $t=\frac{k}{n}$. Therefore,
$$
Q_n\left[\sup_{t\leq\tau_n^q}\left|L_t^n-\int_0^tb(X_s^n)ds\right|\geq\epsilon\right]
\leq\frac{\E_{Q_n}\left[\sup_{t\leq\tau_n^q}\left|L_t^n-\int_0^tb(X_s^n)ds\right|\right]}{\epsilon}.
$$
We can easily prove that $\displaystyle{\lim_{n\rightarrow\infty}\hat{P}\left[\sup_{t\leq\tau_n^q}\left|
L_t^n-\int_0^tb(X_s^n)ds\right|\geq\epsilon\right]= 0}$ as $n$ tends to infinity, as desired. \qed

{\bf Proof of Lemma \ref{l:4.3}.}
 Applying Taylor's expansion to the function
$x^{2m}$, we obtain \beqnn
&&[S_{k+1}^n]^{2m}\\
&=&[S_k^n]^{2m}+2m[S_k^n]^{2m-1}(S_{k+1}^n-S_k^n)+m(2m-1)[\bar{S}_k^n]^{2m-2}(S_{k+1}^n-S_k^n)^2\\
&=&[S_k^n]^{2m}+2m[S_k^n]^{2m-1}\left(\frac{b+\lambda
S}{n}+\frac{\sigma
S}{\sqrt{n}}\varepsilon^{k+1}\right)+m(2m-1)[\bar{S}_k^n]^{2m-2}\left(\frac{b+\lambda
S}{n}+\frac{\sigma S}{\sqrt{n}}\varepsilon^{k+1}\right)^2, \eeqnn
where $$\bar{S}_k^n=S_k^n+\beta\left(\frac{b(
S_k^n)}{n}+\frac{\sigma(S_k^n)}{\sqrt{n}}\varepsilon^{k+1}\right),~~\beta\in[0,1].$$ Moreover, \beqnn \E_{Q_n}[\varepsilon^{k+1}~|~\FF_k^n]
&=&\E_{P_n}\left[\left.\frac{2\tilde{\pi}_{k+1}^n\beta_{k+1}^n}{\tilde{\pi}_k^n\beta_k^n}\varepsilon^{k+1}~\right|~\FF_k^n\right]\\
&=&\E_{P_n}\left[\left.\left(1+\frac{\theta(S_k^n)}{\sqrt{n}}\varepsilon^{k+1}\right)\varepsilon^{k+1}~\right|~\FF_k^n\right]\\
&=&\frac{1}{2}\left(1+\frac{\theta(S_k^n)}{\sqrt{n}}\right)-\frac{1}{2}\left(1-\frac{\theta(S_k^n)}{\sqrt{n}}\right)
=\frac{\theta(S_k^n)}{\sqrt{n}}. \eeqnn Then $$ \frac{b+\lambda
S}{n}+\frac{\sigma
S}{\sqrt{n}}\E_{P_n}\left[\varepsilon^{k+1}~|~\FF_k^n\right]=\frac{b+\lambda
S}{n}+\frac{\sigma S\theta}{n}=\frac{(r+\lambda)S}{n}.$$ Notice that
$$|\bar{S}_k^n|\leq|S_k^n|+|b|+|\lambda S|+|\sigma S|,~|\varepsilon^{k+1}|=1,$$ and $$x^{2m-2}\leq1+x^{2m},~(x+y)^{m}\leq2^m(x^{m}+y^m),~x^2r(x)\leq K(1+x^2)$$ when $x,y>0$.
Taking the conditional expectation with respect to $\FF_k^n$ under $Q_n$,
we have
 \beqnn
\E_{Q_n}\left[[S_{k+1}^n]^{2m}~|~\FF_k^n\right]
&\leq&[S_k^n]^{2m}+2m[S_k^n]^{2m-1}\left(\frac{|r+\lambda
||S|}{n}\right)\\
&&+\frac{m(2m-1)}{n}\left(|S_k^n|+|b|+|\lambda S|+|\sigma
S|\right)^{2m-2}(|b|+|\lambda S|+|\sigma S|)^2.
 \eeqnn
Given the conditions on $b$, $\lambda S$, $\sigma S$, we can find a
constant $K'>0$, such that for any $x\in R$,
 $$|b(x)|\leq K'(1+|x|),~|\lambda(x)x|\leq K'(1+|x|),~|\sigma(x)x|\leq K'(1+|x|),~|b(x)|^2\leq K'(1+x^2),$$
  $$|\lambda(x)x|^2\leq K'(1+x^2),~|\sigma(x)x|^2\leq K'(1+x^2),~\text{and}~|x^2r(x)|\leq
  K'(1+x^2).$$ Hence we can obtain
 \beqnn
&&\E_{Q_n}\left[[S_{k+1}^n]^{2m}~|~\FF_k^n\right]\\
&\leq&[S_k^n]^{2m}+\frac{2mK'}{n}(1+2[S_k^n]^{2m})+\frac{9K'^2m(2m-1)}{n}(3K'+(1+3K')[S_k^n])^{2m-2}(1+[S_k^n]^{2})\\
&\leq&[S_k^n]^{2m}+\frac{2mK'}{n}(1+2[S_k^n]^{2m})+\frac{9(2+6K')^{2m-2}m(2m-1)}{n}(1+[S_k^n]^{2m-2})(1+[S_k^n]^{2})\\
&\leq&K/n+(1+K/n)(1+[S_k^n]^{2m}),
 \eeqnn
and furthermore, $\displaystyle{\E_{Q_n}[[S_k^n]^{2m}~|~\FF_l^n]\leq(1+K/n)^{k-l}(1+[S_l^n]^{2m})\leq
A(1+[S_l^n]^{2m})}$, where $K$ depends on $K'$ and $m$, $0\leq l\leq k$, and
$A=\sup\limits_n(1+K/n)^{k-l}$.

Since $\displaystyle{\E_{Q_n}\left[S_{k+1}^n~|~\FF_k^n\right]\geq\E_{Q_n}\left[\left.S_{k+1}^n\frac{B_k^n}{B_{k+1}^n}~\right|~\FF_k^n\right]=S_k^n}$,
we can get that $(S_k^n)$ is a submartingale. Moreover, $\varphi(x)=x^{2m}$ is a
convex and increasing function in $\IR^+$ and $(S_k^n)$ is
nonnegative. By Jensen's inequality, we have $\displaystyle{\E_{Q_n}[(S_{k+1}^n)^{2m}~|~\FF_k^n]\geq
\E_{Q_n}[S_{k+1}^n~|~\FF_k^n]^{2m}\geq(S_k^n)^{2m}}$. It is easy to see
that $((S_k^n)^{2m})$ is a submartingale. By submartingale
inequality, we have \begin{align*} \E_{Q_n}[[S_k^n]^{2m}~|~\FF_l^n]&\leq
A(1+\E_{Q_n}[\sup_{0\leq t \leq
T}|\tilde{S}_t^n|^{2m}])\leq A(1+\left(\frac{2m}{2m-1}\right)^{2m}\E_{Q_n}[(\tilde{S}_T^n)^{2m}])\\
&\leq A(1+\left(\frac{2m}{2m-1}\right)^{2m}A(1+(S_0^n)^{2m})])\leq C(1+(S_0)^{2m}),
 \end{align*}
where $C$ is large enough, and $\tilde{S}_T^n=S_n^n$.\qed

{\bf Proof of Lemma \ref{l:4.4}.} By Lemma \ref{l:4.1}, $\tilde{Z}^n$ converges
weakly to $Z$, $Y$ is a continuous function of $\tilde{Z}^n$.
Applying continuous mapping theorem, we get \beqlb\label{e:4.2}
Y\left(\tilde{S}_{\cdot}^n,\frac{[n\cdot]}{n}\right)\Rightarrow
Y(S_{\cdot},\cdot),\quad n\rightarrow\infty.\eeqlb Since $\displaystyle{Y_n\left(\tilde{S}_t^n,\frac{[nt]}{n}\right)=Y\left(\tilde{S}_t^n,\frac{[nt]}{n}\right)-\tilde{e}_t^n}$,
where $\displaystyle{\tilde{e}_t^n=Y\left(\tilde{S}_t^n,
\frac{[nt]}{n}\right)-Y_n\left(\tilde{S}_t^n,\frac{[nt]}{n}\right)}$, we
need only prove that the stochastic process $\tilde{e}_\cdot^n$
converges weakly to zero.

Let ``$+$'' and ``$-$'' denote the states ${\varepsilon^{k+1}=1}$ and
${\varepsilon^{k+1}=-1}$ respectively, and define
$S_{k+1}^{n+}=S_{k+1}^n(\omega_1)$,
$S_{k+1}^{n-}=S_{k+1}^n(\omega_2)$. We define two functions as
follows.$$ f_+^{k,n}(t)=Y(S_k^n+t(S_{k+1}^{n+}-S_k^n),
t_k^n+t(t_{k+1}^n-t_k^n)),~~~f_-^{k,n}(t)=Y(S_k^n+t(S_{k+1}^{n-}-S_k^n),
t_k^n+t(t_{k+1}^n-t_k^n)),
$$

Let
$$\frac{\partial Y}{\partial S}=Y_S,~~~\frac{\partial Y}{\partial
t}=Y_t,~~~\frac{\partial^2 Y}{\partial S^2}=Y_{SS},~~~\frac{\partial^2 Y}{\partial t^2}=Y_{tt},~~~\frac{\partial^2
Y}{\partial S\partial t}=Y_{St},~~~Y\left(S_k^n, \frac{k}{n}\right)=Y_k.$$ Then
by Taylor's expansion,
\beqnn f_+(1)&=&f_+(0)+f_+^{'}(0)+\frac{1}{2}f_+^{''}(0)+R_k^n\\
&=&Y_k+Y_S(S_{k+1}^{n+}-S_k^n)+\frac{1}{n}Y_t+\frac{1}{2}Y_{SS}(S_{k+1}^{n+}-S_k^n)^2+\frac{1}
{2n^2}Y_{tt}+Y_{St}\frac{1}{n}(S_{k+1}^{n+}-S_k^n)+R_k^n, \eeqnn
where $\displaystyle{R_k^n=\frac{1}{2}\int_0^1(1-s)^2f_+^{(3)}(s)ds}$. The
expression of $f_-(1)$ is similar to $f_+(1)$ with $S_{k+1}^{n+}$
replaced by $S_{k+1}^{n-}$. By denoting the remaining terms by
$Q_k^n$, we have \beqnn
&&\tilde{\pi}(+;S_k^n)f_+(1)+\tilde{\pi}(-;S_k^n)f_-(1)\\
&=&\left(1+\frac{\tilde{r}_k}{n}\right)^{-1}\left[Y_k+\frac{(r_k+\lambda_k)S_kY_S+Y_t}{n}
+\left(\frac{(b_k+\lambda_kS_k)^2+2(b_k+\lambda_kS_k)\sigma_k
\theta_kS_k}{2n^2}+\frac{\sigma_k^2S_k^2}{2n}\right)Y_{SS}\right.\\
&&\left.+\frac{1}{2n^2}Y_{tt}+\frac{(r_k+\lambda_k)S_k}{n^2}Y_{St}\right]-\gamma_k^n.\eeqnn
By Lemma \ref{l:2.1}, the above equation equals to
\beqnn
&&\left(1+\frac{\tilde{r}_k}{n}\right)^{-1}\left[\left(1+\frac{r_k+\lambda_k}{n}\right)Y_k+
\left(\frac{(b_k+\lambda_kS_k)^2+2(b_k+\lambda_kS_k)\sigma_k
\theta_kS_k}{2n^2}+\frac{Y_{tt}}{2n^2}+\frac{(r_k+ \lambda_k)S_k^n}{n^2}Y_{St}\right)\right]
-\gamma_k^n \\
&=&Y_k-\left(1+\frac{\tilde{r}_k}{n}\right)^{-1}\frac{1}{n^2}m\left(S_k^n,\frac{k}{n}\right)-\gamma_k^n,
\eeqnn where $\gamma_k^n=-\tilde{\pi}(+;S_k^n)R_k^n-\tilde{\pi}(-;S_k^n)Q_k^n$, and
$$m\left(S_k^n,\frac{k}{n}\right)=\left(\frac{(b_k+\lambda_kS_k)^2+2(b_k+\lambda_kS_k)\sigma_k
\theta_kS_k}{2n^2}+\frac{1}{2n^2}Y_{tt}+\frac{(r_k+\lambda_k)S_k^n}{n^2}Y_{St}\right)+Y_k\left(\frac{\lambda_k^2+2r_k\lambda_k}{2}+
o\left(\frac{1}{n}\right)\right).$$ Hence we obtain the following
recurrent equation for $e_k^n$,
 \beqnn
e_k^n&=&\tilde{\pi}(+;S_k^n)Y\left(S_{k+1}^{n+},\frac{k+1}{n}\right)+\tilde{\pi}(-;S_k^n)Y\left(S_{k+1}^{n-},\frac{k+1}{n}\right)
+\left(1+\frac{\tilde{r}_k}{n}\right)^{-1}\frac{1}{n^2}m\left(S_k^n,\frac{k}{n}\right)
\\
&&+\gamma_k^n - \tilde{\pi}(+;S_k^n)Y_n\left(S_{k+1}^{n+},
\frac{k+1}{n}\right)-\tilde{\pi}(-;S_k^n)Y_n\left(S_{k+1}^{n+},\frac{k+1}{n}\right)\\
&=&\tilde{\pi}(+;S_k^n)e_{k+1}^{n+}+\tilde{\pi}(-;S_k^n)e_{k+1}^{n-}+
\left(1+\frac{\tilde{r}_k}{n}\right)^{-1}\frac{1}{n^2}m\left(S_k^n,\frac{k}{n}\right)+\gamma_k^n.
\eeqnn By the definition of $\tilde{\pi}(\cdot;S_k^n)$, we obtain
$$
e_k^n=\E_{Q_n}\left[\left.e_{k+1}^n\left(1+\frac{\tilde{r}_k}{n}\right)^{-1}~\right|~\FF_k^n\right]+
\left(1+\frac{\tilde{r}_k}{n}\right)^{-1}\frac{1}{n^2}m\left(S_k^n,\frac{k}{n}\right)+\gamma_k^n.
$$
Since $e_n^n=Y(S_n^n,1)-Y_n(S_n^n,1)=g(S_n^n)-g(S_n^n)=0$, we get
 \beqnn
e_k^n=\E_{Q_n}\left[\left.\sum_{i=k}^{n-1}\frac{1}{n^2}m\left(S_i^n,\frac{i}{n}\right)
\frac{\beta_k^n}{\beta_{i+1}^n}+
\gamma_i^n\frac{\beta_k^n}{\beta_i^n}~\right|~\FF_k^n\right].
 \eeqnn

By the assumption that $Y$ and its derivative satisfy linear
increasing condition, there exists constants $C_1>0$ and $q$, such
that $\displaystyle{\left|m\left(S_i^n,\frac{i}{n}\right)\right|\leq C_1(1+|S_i^n|^{2q})}$.

By Lemma \ref{l:4.3}, for $k\leq i\leq n$ there exists constant $C>0$
large enough such that $\displaystyle{E_{Q_n}[|S_i^n|^{2q}~|~\FF_k^n]\leq
C(1+|S_0|^{2q})}$.
Therefore,
\begin{align*}
&\E_{Q_n}\left[\left.\sum_{i=k}^{n-1}\left|\frac{1}{n^2}m\left(S_i^n,\frac{i}{n}\right)
\frac{\beta_k^n}{\beta_{i+1}^n}\right|~\right|~\FF_k^n\right]\leq\frac{1}{n^2}\sum_{i=k}^{n-1}\E_{Q_n}\left[\left.\left|m\left(S_i^n,\frac{i}{n}\right)\right|~\right|~\FF_k^n\right]\\
\leq&\frac{1}{n^2}\sum_{i=k}^{n-1}\E_{Q_n}[C_1(1+|S_i^n|^{2q})~|~\FF_k^n]\leq\frac{1}{n^2}\sum_{i=k}^{n-1}C_1(1+C(1+|S_0|^{2q})).
\end{align*}

For the second part we can also write out the expressions of
$f_+^{(3)}(s)$, $f_-^{(3)}(s)$, they are of order
$n^{-\frac{3}{2}}$, by analogous argument we can choose $q$ large
enough and constant $D>0$ satisfying $$
\E_{Q_n}\left[\left.\sum_{i=k}^{n-1}\left|\gamma_i^n\frac{\beta_k^n}{\beta_i^n}\right|~\right|~S_k^n\right]
\leq\frac{D}{\sqrt{n}}(1+|S_k^n|^{2q}). $$ Then we can choose
$\tilde{C}$ large enough which depends on $q$, $k$ and $n$, such that $\displaystyle{|e_k^n| \leq
\frac{\tilde{C}}{\sqrt{n}}(1+|S_0|^{2q})}$. So $$
\hat{P}\left(\sup_{0\leq t \leq 1}|\tilde{e}_t^n|\geq
\epsilon\right)\leq\hat{P}\left(\sup_{0\leq t \leq
1}\frac{\tilde{C}(1+|\tilde{S}_0|^{2q})}{\sqrt{n}}\geq\epsilon\right)\leq\frac{\tilde{C}\left(1+\sup\limits_{0\leq t \leq
1}|\tilde{S}_0|^{2q}\right)}{\epsilon\sqrt{n}}\rightarrow0$$ as $n\rightarrow\infty$. Therefore,
$\sup\limits_{0\leq t \leq 1}|\tilde{e}_t^n|\rightarrow 0$ as
$n\rightarrow\infty$, which means that $\tilde{e}_t^n$ converges to zero almost surely, that is,
$\tilde{e}_{\cdot}^n\Rightarrow0$. Combined with (\ref{e:4.2}),
we get the conclusion. \qed

\section{Conclusion}

In this paper, the weak convergence of discrete--time equity
derivatives pricing model with default risk is proved in a
no--arbitrage framework. Our results present a mathematical
foundation for derivative pricing with default risk using numerical
method. It remains to study the convergence for the hedging
strategy.

\bigskip
\noindent{\bf Acknowledgements.}  We would like to thank the anonymous referee for the useful suggestions, which are a great help to improve the manuscript. We thank Xiaobin Li for his careful reading of the manuscript and useful comments. The first author would like to thank
Professors Zenghu Li and Wenming Hong for their invaluable comments and
encouragements during her stay in Beijing Normal University.


\bigskip


\begin{thebibliography}{99}

\bibitem{aldous1978stopping}  D. Aldous, {\it Stopping times and tightness},
Ann. Probab, 10 (1978), 335--340.

\bibitem{bielecki2008defaultable}  T. Bielecki, S. Cr$\acute{e}$pey, M. Jeanblanc and
M. Rutkowski, {\it Defaultable options in a Markovian intensity
model of credit risk}, Math. Finance, 18 (2009), 493-518.

\bibitem{bielecki2004credit} T. Bielecki and M. Rutkowski, {\it Credit
Risk: Modeling, Valuation and Hedging}, Springer-Verlag, Berlin,
Heidelberg, New York, 2002.


\bibitem{blanchet2004hazard} C. Blanchet-Scalliet and M. Jeanblanc, {\it Hazard
rate for credit risk and hedging defaultable contingent claims},
Finance Stoch., 8 (2004), 145-159.

%

\bibitem{cox1979option} J. Cox and M. Rubinstein, Option pricing: A
simplied approach, J. Finan. Econ., 7 (1979), 229-263.


\bibitem{duffie2006discrete} D. Duffie and P. Protter, {\it From Discrete to Continuous Time
Finance: Weak convergence of the Financial Gain Process},
Mathematical Finance, 2 (1992), 1-15.

\bibitem{duffie1999modeling}  D. Duffie and K. Singleton, {\it Modeling term
structure of defaultable bonds}, Revie. Finan. Studies, 12 (1999),
687-720.

\bibitem{ethier1986markov} S. N. Ethier and T. G. and Kurtz, {\it Markov
Processes: Characterization and Convergence}, John Wiley, New York,
1986.

\bibitem{he1990convergence} H. He, {\it Convergence from discrete to continuous
time contingent claims prices}, Revie. Finan. Studies,  3 (1989),
523-546.

\bibitem{jacod1987limit} J. Jacod,  and A. N. Shiryaev, {\it Limit
Theorems for Stochastic Processes}, Springer-Verlag, Berlin,
Heidelberg, New York, 1987.

\bibitem{jarrow2012pricing} R. Jarrow and S. Turnbull,  {\it Pricing Options
on Financial Securities  Subject  to Credit Risk},  J. Finance, 50
(1995), 53--85.
%
%
%

\bibitem{lando1998cox}  D. Lando, {\it On Cox Processes and Credit Risky
Securities}, Revie. Deriv. Research, 2 (1998), 99-120.

%

\bibitem{nieuwenhuis2004weak} J. W. Nieuwenhuis and M. H. Vellekoop, {\it Weak convergence of
three methods, to price options on defautable assets}, Decis.
Econom. Finance, 27, (2004), 87-107.


\end{thebibliography}
\end{document}